\documentclass[11pt]{article}
\usepackage{graphicx}
\usepackage{color}
\usepackage{amsmath}
\usepackage{amssymb}
\usepackage{amscd}
\usepackage{bbm}
\newcommand{\R}{\mathbb{R}}
\newcommand{\inr}[1]{\bigl< #1 \bigr>}

\newcommand{\E}{\mathbb{E}}

\newcommand{\eps}{\varepsilon}

\newtheorem{Theorem}{Theorem}[section]
\newtheorem{Lemma}[Theorem]{Lemma}

\newtheorem{Definition}[Theorem]{Definition}

\newtheorem{Remark}[Theorem]{Remark}

\newtheorem{Question}[Theorem]{Question}
\numberwithin{equation}{section}

\newcommand{\norm}[1]{\left\|#1\right\|}%

\def \proof {\noindent {\bf Proof.}\ \ }

\def \endproof
{{\mbox{}\nolinebreak\hfill\rule{2mm}{2mm}\par\medbreak}}
\def\IND{\mathbbm{1}}

\begin{document}
\title{{On multiplier processes under weak moment assumptions}}
\author{Shahar Mendelson${}^{1,2}$}

\footnotetext[1]{Department of Mathematics, Technion, I.I.T., Haifa, Israel and Mathematical Sciences Institute, The Australian National University, Canberra, Australia, Email:
shahar@tx.technion.ac.il}
\footnotetext[2]{Supported in part by the Israel Science Foundation.}

\maketitle

\begin{abstract}
We show that if $V \subset \R^n$ satisfies a certain symmetry condition (closely related to unconditionaity) and if $X$ is an isotropic random vector for which $\|\inr{X,t}\|_{L_p} \leq L \sqrt{p}$ for every $t \in S^{n-1}$ and $p \lesssim \log n$, then the corresponding empirical and multiplier processes indexed by $V$ behave as if $X$ were $L$-subgaussian.
\end{abstract}

\section{Introduction}
The motivation for this work comes from various problems in Learning Theory, in which one encounters the following random process.

Let $X=(x_1,...,x_n)$ be a random vector on $\R^n$ (whose coordinates $(x_i)_{i=1}^n$ {\it need not} be independent) and let $\xi$ be a random variable that need not be independent of $X$. Set $(X_i,\xi_i)_{i=1}^N$ to be $N$ independent copies of $(X,\xi)$, and for $V \subset \R^n$ define the centred multiplier process
\begin{equation} \label{eq:multi-def-intorduction}
\sup_{v \in V} \left|\frac{1}{\sqrt{N}} \sum_{i=1}^N (\xi_i \inr{X_i,v} - \E \xi \inr{X,v}) \right|.
\end{equation}
Multiplier processes are often studied in a more general context, in which the indexing class need not be a class of linear functionals on $\R^n$. Instead, one may consider an arbitrary probability space $(\Omega,\mu)$ and in which case $F$ is a class of functions on $\Omega$. Let $X_1,...,X_N$ be independent, distributed according to $\mu$, and the multiplier process indexed by $F$ is
\begin{equation} \label{eq:multi-general-def-intorduction}
\sup_{f \in F} \left|\frac{1}{\sqrt{N}} \sum_{i=1}^N (\xi_i f(X_i) - \E \xi f(X_i)) \right|.
\end{equation}
Naturally, the simplest multiplier process is when $\xi \equiv 1$ and \eqref{eq:multi-general-def-intorduction} is the standard empirical process.

\vskip0.4cm

Controlling a multiplier process is relatively straightforward when $\xi \in L_2$ and is independent of $X$. For example, one may show (see, e.g., \cite{vanderVaartWellner}, Chapter 2.9) that if $\xi$ is a mean-zero random variable that is independent of $X_1,...,X_N$ then
$$
\E \sup_{f \in F} \left|\frac{1}{\sqrt{N}} \sum_{i=1}^N (\xi_i f(X_i) - \E \xi f(X_i)) \right| \leq C \|\xi\|_{L_2} \E \sup_{f \in F} \left|\frac{1}{\sqrt{N}}\sum_{i=1}^N \eps_i f(X_i) \right|,
$$
where here and throughout the article, $(\eps_i)_{i=1}^N$ are independent, symmetric $\{-1,1\}$-valued random variables that are independent of $(X_i,\xi_i)_{i=1}^N$, and $C$ is an absolute constant.

This estimate and others of its kind show that multiplier processes are as `complex' as their seemingly simpler empirical counterparts. However, the results we are looking for are of a different nature: estimates on multiplier processes that are based on some natural complexity parameter of the underlying class $F$, and that exhibits the class' geometry.

It turns out that chaining methods lead to such estimates, and the structure of $F$ may be captured using the following parameter, which is a close relative of Talagrand's $\gamma$-functionals \cite{MR3184689}.
\begin{Definition} \label{def:Lambda}
For a random variable $Z$ and $p \geq 1$, set
$$
\|Z\|_{(p)} = \sup_{1 \leq q \leq p} \frac{\|Z\|_{L_q}}{\sqrt{q}}.
$$
Given a class of functions $F$, $u \geq 1$ and $s_0 \geq 0$, put
\begin{equation} \label{eq:comp-gamma-s-0}
{\Lambda}_{s_0,u}(F) = \inf \sup_{f \in F} \sum_{s \geq s_0} 2^{s/2} \|f - \pi_s f\|_{(u^22^s)},
\end{equation}
where the infimum is taken with respect to all sequences $(F_s)_{s \geq 0}$ of subsets of $F$, and of cardinality $|F_s| \leq 2^{2^s}$. $\pi_s f$ is the nearest point in $F_s$ to $f$ with respect to the $(u^22^s)$ norm.

Let
$$
\tilde{\Lambda}_{s_0,u}(F)={\Lambda}_{s_0,u}(F)+2^{s_0/2}\sup_{f \in F} \|\pi_{s_0}f\|_{(u^22^{s_0})}.
$$
\end{Definition}
To put these definitions in some perspective, $\|Z\|_{(p)}$ measures the local-subgaussian behaviour of $Z$, and the meaning of `local' is that $\| \ \|_{(p)}$ takes into account the growth of $Z$'s moments up to a fixed level $p$. In comparison,
$$
\|Z\|_{\psi_2} \sim \sup_{q \geq 2} \frac{\|Z\|_{L_q}}{\sqrt{q}},
$$
implying that for $2 \leq p <\infty$, $\|Z\|_{(p)} \lesssim \|Z\|_{\psi_2}$; hence, for every $u \geq 1$ and $s \geq s_0$,
$$
{\Lambda}_{s_0,u}(F) \lesssim \inf \sup_{f \in F} \sum_{s \geq s_0} 2^{s/2} \|f - \pi_s f\|_{\psi_2},
$$
and $\tilde{\Lambda}_{0,u}(F) \leq c\gamma_2(F,\psi_2)$ (see \cite{MR3184689} for a detailed study on {\it generic chaining} and the $\gamma$ functionals).

Recall that the canonical gaussian process indexed by $F$ consists of centred gaussian random variable $G_f$, and the covariance structure of the process is endowed by the inner product in $L_2(\mu)$. Let
$$
\E \sup_{f \in F} G_f= \sup \{ \E \sup_{f  \in F^\prime} G_f : \ F^\prime \subset F, \ F^\prime \ {\rm is \ finite} \}.
$$
and note that if the class $F \subset L_2(\mu)$ is $L$-subgaussian, that is, if for every $f,h \in F \cup 0$, $$
\|f-h\|_{\psi_2(\mu)} \leq L \|f-h\|_{L_2(\mu)},
$$
then $\tilde{\Lambda}_{s_0,u}(F)$ may be bounded using the canonical gaussian process indexed by $F$. Indeed, by Talagrand's Majorizing Measures Theorem \cite{MR906527,MR3184689}, for every $s_0 \geq 0$,
$$
\tilde{\Lambda}_{s_0,u}(F) \lesssim L\bigl(\E \sup_{f \in F} G_f + 2^{s_0/2}\sup_{f \in F} \|f\|_{L_2(\mu)}\bigr).
$$
As an example, let $V \subset \R^n$ and set $F=\{\inr{v,\cdot} : v \in V\}$ to be the class of linear functionals endowed by $V$. If $X$ is an isotropic, $L$-subgaussian vector, it follows that for every $t \in \R^n$,
$$
\|\inr{X,t}\|_{\psi_2} \leq L\|\inr{X,t}\|_{L_2} = L \|t\|_{\ell_2^n}.
$$
Therefore, if $G=(g_1,...,g_n)$ is the standard gaussian vector in $\R^n$, $\ell_*(V)=\E \sup_{v \in V} |\inr{G,v}|$ and $d_2(V)=\sup_{v \in V} \|v\|_{\ell_2^n}$, one has
\begin{align*}
\tilde{\Lambda}_{s_0,u}(F) \lesssim & L\bigl(\E \sup_{v \in V} \inr{G,v} + 2^{s_0/2}\sup_{v \in V} \|\inr{X,v}\|_{L_2} \bigr)
\\
\lesssim & L \bigl(\ell_*(V) + 2^{s_0/2}d_2(V)\bigr).
\end{align*}

As the following estimate from \cite{shahar_multi_pro} shows, $\tilde{\Lambda}$ can be used to control a multiplier process in a relatively general situation.

\begin{Theorem} \label{thm:multiplier-intro}
For $q>2$, there are constants $c_0$, $c_1,c_2,c_3$ and $c_4$ that depend only on $q$ for which the following holds. Let $\xi \in L_q$ and set $\xi_1,...,\xi_N$ to be independent copies of $\xi$. Fix an integer $s_0 \geq 0$ and $w,u>c_0$. Then, with probability at least
$$
1-c_1w^{-q} N^{-((q/2)-1)}\log^{q} N-2\exp(-c_2u^2 2^{s_0}),
$$
$$
\sup_{f \in F} \left|\frac{1}{\sqrt{N}} \sum_{i=1}^N \left(\xi_i f(X_i) - \E \xi f\right) \right| \leq c_3wu\|\xi\|_{L_q}\tilde{\Lambda}_{s_0,c_4u}(F).
$$
\end{Theorem}

It follows from Theorem \ref{thm:multiplier-intro} that if
$$
D(V) = \left(\frac{\ell_*(V)}{d_2(V)}\right)^2
$$
then with probability at least
$$
1-c_2w^{-q} N^{-((q/2)-1)}\log^{q} N-2\exp(-c_3u^2 D(V)),
$$
\begin{equation} \label{eq:multi-subgaussian-intro}
\sup_{f \in F} \left|\frac{1}{\sqrt{N}}\sum_{i=1}^N \left(\xi_i \inr{v,X_i} - \E \xi \inr{v,X}\right) \right| \lesssim Lwu\|\xi\|_{L_q}\ell_*(V).
\end{equation}

There are other generic situations in which $\tilde{\Lambda}_{s_0,u}(F)$ may be controlled using the geometry of $F$ (for example \cite{MR2899978,shahar_multi_pro} when $F$ is a class of linear functionals on $\R^n$ and $X$ is an unconditional, log-concave random vector). However, there is no satisfactory theory that describes $\tilde{\Lambda}_{s_0,u}(F)$ for an arbitrary class $F$; such results are highly nontrivial.

Moreover, because the definition of $\Lambda_{s_0,u}(F)$ involves $\| \ \|_{(p)}$ for every $p$, class members must have arbitrarily high moments for $\Lambda_{s_0,u}$ to be well defined.

\vskip0.4cm

In the context of classes of linear functionals on $\R^n$, one expects an analogous result to Theorem \ref{thm:multiplier-intro} to be true even if the functionals $\inr{X,t}$ do not have arbitrarily high moments. A realistic conjecture is that if for each $t \in S^{n-1}$
$$
\|\inr{X,t}\|_{L_q} \leq L\sqrt{q} \|\inr{X,t}\|_{L_2} \ \ {\rm  for \ every} \ 2 \leq q \lesssim n
$$
then a subgaussian-type estimate like \eqref{eq:multi-subgaussian-intro} should still be true.

In what follows we will not focus on such a general result that is likely to hold for {\it every} $V \subset \R^n$. Rather, we will concentrate our attention on situations where a subgaussian estimate like \eqref{eq:multi-subgaussian-intro} is true, but linear functionals only satisfy
$$
\|\inr{X,t}\|_{L_q} \leq L\sqrt{q} \|\inr{X,t}\|_{L_2} \ \ {\rm for \ every} \  2 \leq q \lesssim \log n.
$$
The obvious example in which only $\sim \log n$ moments should suffice is $V=B_1^n$ (or similar sets that have $\sim n$ extreme points). Having said that, the applications that motivated this work require a broader spectrum of sets that only need that number of moments to exhibit a subgaussian behaviour as in \eqref{eq:multi-subgaussian-intro}.

\begin{Question} \label{qu:main}
Let $X=(x_1,...,x_n)$ be an isotropic random vector and assume that $\|x_i\|_{L_q} \leq L\sqrt{q}$ for every $2 \leq q \leq p$. If $\xi \in L_{q_0}$ for some $q_0 >2$, how small can $p$ be while still having that
$$
\E \sup_{v \in V} \left| \frac{1}{\sqrt{N}} \sum_{i=1}^N \xi_i \inr{X_i,v} - \E \xi \inr{X,v} \right| \leq C(L,q_0) \|\xi\|_{L_{q_0}} \ell_*(V)?
$$
\end{Question}
We will show $p \sim \log{n}$ suffices for a positive answer to Question \ref{qu:main} if the norm $\|z\|_{V^\circ}=\sup_{v \in V} |\inr{v,z}|$ satisfies the following unconditionality property:

\begin{Definition} \label{def:K-unconditional}
Given a vector $x=(x_i)_{i=1}^n$, let $(x_i^*)_{i=1}^n$ be the non-increasing rearrangement of $(|x_i|)_{i=1}^n$.

The normed space $(\R^n,\| \ \|)$ is $K$-unconditional with respect to the basis $\{e_1,...,e_n\}$ if for every $x \in \R^n$ and every permutation of $\{1,...,n\}$
$$
\|\sum_{i=1}^n x_i e_i\| \leq K \|\sum_{i=1}^n x_{\pi(i)}e_i\|,
$$
and if $y \in \R^n$ and $x_i^* \leq y_i^*$ for $1 \leq i \leq n$ then
$$
\|\sum_{i=1}^n x_i e_i\| \leq K \|\sum_{i=1}^n y_ie_i\|
$$
\end{Definition}

\begin{Remark}
This is not the standard definition of an unconditional basis, though every unconditional basis (in the classical sense) on an infinite dimensional space  satisfies Definition \ref{def:K-unconditional} for some constant $K$ (see, e.g., \cite{MR2192298}).
\end{Remark}

There are many natural examples of $K$-unconditional norms, including all the $\ell_p$ norms. Moreover,  the norm $\sup_{v \in V} \sum_{i=1}^n v_i^* z_i^*$ is $1$-unconditional. In fact, if $V \subset \R^n$ is closed under permutations and reflections (sign-changes), then $\|\cdot\|_{V^\circ}$ is $1$-unconditional. Finally, since the maximum of two $K$-unconditional norms is $K$-unconditional, it follows that if $\| \cdot \|_{V^\circ}$ is $K$-unconditional, so is the norm $\sup_{v \in V \cap rB_2^n} \inr{\cdot,v}$.

\vskip0.4cm
We will show the following:
\begin{Theorem} \label{thm:main-formulation}
There exists an absolute constant $c_1$ and for $K \geq 1$, $L \geq 1$ and $q_0 >2$ there exists a constant $c_2$ that depends only on $K$, $L$ and $q_0$ for which the following holds. Consider
\begin{description}
\item{$\bullet$} $V \subset \R^n$ for which the norm $\|\cdot\|_{V^\circ}=\sup_{v \in V} |\inr{v,\cdot}|$ is $K$-unconditional with respect to the basis $\{e_1,...,e_n\}$.
\item{$\bullet$} $\xi \in L_{q_0}$ for some $q_0>2$.
\item{$\bullet$} An isotropic random vector $X \in \R^n$ which satisfies that
$$
\max_{1 \leq j \leq n} \|\inr{X,e_j}\|_{(p)} \leq L \ \ {\rm  for}  \ p =c_1 \log n.
$$
\end{description}
If $(X_i,\xi_i)_{i=1}^N$ are independent copies of $(X,\xi)$ then
$$
\E \sup_{v \in V} \left|\frac{1}{\sqrt{N}} \sum_{i=1}^N \left(\xi_i \inr{X_i,v} - \E \xi \inr{X,v}\right) \right| \leq c_2 \|\xi\|_{L_q} \ell_*(V).
$$
\end{Theorem}

The proof of Theorem \ref{thm:main-formulation} is based on the study of a conditioned Bernoulli process. Indeed, a standard symmetrization argument (see, e.g., \cite{LT:91,vanderVaartWellner}) shows that if $(\eps_i)_{i=1}^N$ are independent, symmetric, $\{-1,1\}$-valued random variables that are independent of $(X_i,\xi_i)_{i=1}^N$ then
$$
\E \sup_{v \in V} \left| \frac{1}{\sqrt{N}} \sum_{i=1}^N \xi_i \inr{X_i,v} - \E \xi \inr{X,v} \right| \leq C \E\sup_{v \in V} \left| \frac{1}{\sqrt{N}} \sum_{i=1}^N \eps_i \xi_i \inr{X_i,v} \right|
$$
for an absolute constant $C$; a similar bound hold with high probability, showing that it suffices to study the supremum of the conditioned Bernoulli process
$$
\sup_{v \in V} \left| \frac{1}{\sqrt{N}} \sum_{i=1}^N \eps_i \xi_i \inr{X_i,v} \right| = (*).
$$
Put $x_i(j) = \inr{X_i,e_j}$ and set $Z_j = N^{-1/2}\sum_{i=1}^N \eps_i \xi_i x_i(j)$, which is a sum of iid random variables. Therefore, if $Z=(Z_1,...,Z_n)$ then
$$
(*) = \sup_{v \in V} \inr{Z,v}.
$$
The proof of Theorem \ref{thm:main-formulation} follows by showing that for a well-chosen constant $C(L,q)$ the event
$$
\left\{Z_j^* \leq C\E g_j^* \ {\rm for \ every \ } 1\leq j \leq n\right\}
$$
is of high probability, and if the norm $\|\cdot\|_{V^\circ}=\sup_{v \in V} \inr{\cdot,v}$ is $K$-unconditional then
$$
\sup_{v \in V} \inr{Z,v} \leq C_1(K,L,q) \E \sup_{v \in V} \inr{G,v}.
$$

\vskip0.4cm

Before presenting the proof of Theorem \ref{thm:main-formulation}, let us turn to one of its outcomes -- estimates on the random Gelfand widths of a convex body. We will present another application, motivated by a question in the rapidly developing area of {\it Spare Recovery} in Section \ref{sec:sparse}.

\vskip0.4cm

Let $V \subset \R^n$ be a convex, centrally symmetric set. A well known question in Asymptotic Geometric Analysis has to do with the diameter of a random $m$-codimensional section of $V$ (see, e.g., \cite{MR827766,MR941809,MR845980,MR3331351}). In the past, the focus was on obtaining such estimates for subspaces selected uniformly according to the Haar measure, or alternatively, according to the measure endowed via the kernel of an $m \times n$ gaussian matrix (see, e.g. \cite{MR1036275}). More recently, there has been a growing interest in other notions of randomness, most notably, generated by kernels of other random matrix ensembles. For example, the following was established in \cite{MR2373017}:
\begin{Theorem} \label{thm:MPT}
Let $X_1,...,X_m$ be distributed according to an isotropic, $L$-subgaussian random vector on $\R^n$, set $\Gamma=\sum_{i=1}^m \inr{X_i,\cdot}e_i$ and put
$$
r_G(V,\gamma) = \inf\{ r>0 : \ell_*(V \cap r B_2^n) \leq \gamma r \sqrt{m} \}.
$$
Then, with probability at least $1-2\exp(-c_1(L)m)$
$$
{\rm diam}({\rm ker}(\Gamma) \cap V) \leq r_G(V,c_2(L)),
$$
for constants $c_1$ and $c_2$ that depends only on $L$.
\end{Theorem}

A version of Theorem \ref{thm:MPT} was obtained under a much weaker assumption: the random vector need not be $L$-subgaussian; rather, it suffices that it satisfies a weak small-ball condition.

\begin{Definition} \label{def:small-ball}
The isotropic random vector $X$ satisfies a small-ball condition with constants $\kappa>0$ and $0<\eps \leq 1$ if for every $t \in S^{n-1}$,
$$
Pr(|\inr{X,t}| \geq \kappa) \geq \eps.
$$
\end{Definition}
The analog of gaussian parameter $r_G$ for a general random vector $X$ is
$$
r_X(V,\gamma) = \inf\Bigl\{ r>0 : \E \sup_{v \in V \cap rB_2^n} \bigl| \frac{1}{\sqrt{m}}\sum_{i=1}^m \inr{X_i,v} \bigr| \leq \gamma r \sqrt{m} \Bigr\}.
$$
Clearly, if $X$ is $L$-subgaussian then $r_X(V,\gamma) \leq r_G(V,cL \gamma)$ for a suitable absolute constant $c$.
\begin{Theorem} \label{thm:heavy-tailed-gelfand} \cite{MR3364699,shahar_general_loss}
Let $X$ be an isotropic random vector that satisfies the small-ball condition with constants $\kappa$ and $\eps$. If $X_1,...X_m$ are independent copies of $X$ and $\Gamma=\sum_{i=1}^m \inr{X_i,\cdot}e_i$, then
with probability at least $1-2\exp(-c_0(\eps) m)$
$$
{\rm diam}({\rm ker}(\Gamma) \cap V) \leq r_X\bigl(V,c_1(\kappa,\eps)\bigr).
$$
\end{Theorem}

Theorem \ref{thm:main-formulation} implies that if the norm $\|z\|_{V^\circ}$ is $K$-unconditional, and the growth of moments of the coordinate linear functionals $\inr{X,e_i}$ for $1 \leq i \leq n$ is $L$-`subgaussian' up to the level $\sim \log n$, then the small-ball condition depends only on $L$ and $r_X(V,c_1(L)) \leq r_G(V,c_2(L,K))$. Therefore, with probability at least $1-2\exp(-c_0(L) m)$ one has the gaussian estimate:
$$
{\rm diam}({\rm ker}(\Gamma) \cap V) \leq r_G\bigl(V,c_2(L,K)\bigr),
$$
even though the choice of a subspace has been made according to an ensemble that could be very far from a subgaussian one.

\vskip0.4cm
We end this introduction with a word about notation. Throughout, absolute constants are denoted by $c,c_1...$, etc. Their value may change from line to line or even within the same line. When a constant depends on a parameter $\alpha$ it will be denoted by $c(\alpha)$. $A \lesssim B$ means that $A \leq cB$ for an absolute constant $c$, and the analogous two-sided inequality is denoted by $A \sim B$. In a similar fashion, $A \lesssim_\alpha B$ implies that $A \leq c(\alpha)B$, etc.

\section{Proof of Theorem \ref{thm:main-formulation}}
There are two substantial difficulties in the proof of Theorem \ref{thm:main-formulation}. First, $Z_1,...,Z_n$ are not independent random variables, not only because of the Bernoulli random variables $(\eps_i)_{i=1}^N$ that appear in all the $Z_i$'s, but also because the coordinates of $X=(x_1,...,x_n)$ need not be independent. Second, while there is some flexibility in the moment assumptions on the coordinates of $X$, there is no flexibility in the moment assumption on $\xi$, which is only `slightly better' than square-integrable.

\vskip0.4cm

As a starting point, let us address the fact that the coordinates of $Z$ need not be independent.

\begin{Lemma} \label{lemma:W-j-union-bound}
There exist absolute constants $c_1$ and $c_2$ for which the following holds. Let $\beta \geq 1$ and set $p=2\beta \log(en)$. If $(W_j)_{j=1}^n$ are random variables and satisfy that $\|W_j\|_{(p)} \leq L$, then for every $t \geq 1$, with probability at least $1-c_1t^{-2\beta}$,
$$
W_j^* \leq c_2 t  L \sqrt{\beta \log(en/j)} \ \ \ {\rm for \ every \ } 1 \leq j \leq n.
$$
\end{Lemma}

\proof Let $a_1,...,a_k \in \R$ and by the convexity of $t \to t^q$,
$$
\bigl(\frac{1}{k}\sum_{j=1}^k a_j^2\bigr)^q \leq \frac{1}{k} \sum_{j=1}^k a_j^{2q}.
$$
Thus, given $(a_i)_{i=1}^n$, and taking the maximum over subsets of $\{1,...,n\}$ of cardinality $k$,
$$
\max_{|J_1| = k}\bigl(\frac{1}{k}\sum_{j \in J_1} a_j^2\bigr)^q \leq \max_{|J_1| = k} \frac{1}{k} \sum_{j \in J_1} a_j^{2q} \leq \frac{1}{k} \sum_{j=1}^n a_j^{2q}.
$$
When applied to $a_j = W_j$, it follows that point-wise,
\begin{equation} \label{eq-in-covering-estimate}
\bigl(\frac{1}{k}\sum_{j=1}^k (W_j^*)^2\bigr)^q \leq \frac{1}{k}\sum_{i=1}^n W_j^{2q}.
\end{equation}
Since $\|W_j\|_{(p)} \leq L$ it is evident that $\E W_j^{2q} \leq L^{2q} q^q$ for $2q \leq p$. Hence, taking the expectation in \eqref{eq-in-covering-estimate},
$$
\Bigl(\E \bigl(\frac{1}{k}\sum_{j=1}^k (W_j^*)^2\bigr)^q \Bigr)^{1/q} \leq q L^{2} \cdot \bigl(\frac{n}{k}\bigr)^{1/q} \leq c_1qL^2
$$
for $q = \beta \log (e n/k)$ (which does satisfy $2q \leq p$). Hence, by Chebyshev's inequality, for $t \geq 1$,

\begin{equation} \label{eq:2-in-proof}
Pr\Bigl(\frac{1}{k}\sum_{j \leq k} (W_j^*)^2 \geq  (et)^2 c_1^2 L^2 q \Bigr) \leq \frac{1}{t^{2q}} \cdot e^{-2q} = \left(\frac{k}{en}\right)^2 \cdot \frac{1}{t^{-2q}}.
\end{equation}
Using \eqref{eq:2-in-proof} for $k=2^j$ and applying the union bound, it is evident that with probability at least $1-2t^{-2\beta}$, for every $1 \leq k \leq n$,
$$
(W_k^*)^2 \leq \frac{1}{k}\sum_{j \leq k} (W_j^*)^2 \lesssim  t^2 L^2 \beta \log(en/k).
$$
\endproof

Recall that $q_0>2$ and set $\eta=(q_0-2)/4$. Let $u \geq 2$ and consider the event
$$
{\cal A}_u=\{\xi_i^* \leq u\|\xi\|_{L_{q_0}}(eN/i)^{1/q_0} \ {\rm for \ every \ } 1 \leq i \leq N \}.
$$
A standard binomial estimate combined with Chebyshev's inequality for $|\xi|^{q_0}$ shows that ${\cal A}_u$ is a nontrivial event. Indeed,
\begin{equation*}
Pr\left(\xi_i^* \geq u\|\xi\|_{L_{q_0}}(eN/i)^{1/q_0}\right) \leq \binom{N}{i} Pr^i\left(\xi \geq u\|\xi\|_{L_{q_0}}(eN/i)^{1/q_0}\right) \leq \frac{1}{u^{iq_0}},
\end{equation*}
and by the union bound for $1 \leq i \leq n$, $Pr({\cal A}_u) \leq 2/u^{q_0}$.

\vskip0.4cm

The random variables we shall use in Lemma \ref{lemma:W-j-union-bound} are
$$
W_j=Z_j\IND_{{\cal A}_u},
$$
for $u \geq 2$ and $1 \leq j \leq n$.

\vskip0.4cm

The following lemma is the crucial step in the proof of Theorem \ref{thm:main-formulation}.
\begin{Lemma} \label{lemma:(p)-estimates}
There exists an absolute constant $c$ for which the following holds.
Let $X$ be a random variable that satisfies $\|X\|_{(p)} \leq L$ for some $p>2$ and set $X_1,...,X_N$ to be independent copies if $X$.
If
$$
W=\left|\frac{1}{\sqrt{N}} \sum_{i=1}^N \eps_i \xi_i X_i\right| \IND_{{\cal A}_u},
$$
then $\|W\|_{(p)} \leq c u L$.
\end{Lemma}

The proof of Lemma \ref{lemma:(p)-estimates} requires two preliminary estimates on the `gaussian' behaviour of a monotone rearrangements of $N$ copies of a random variable.

\begin{Lemma} \label{Lemma:com-with-gaussians}
There exists an absolute constant $c$ for which the following holds. Assume that $\|X\|_{(2p)} \leq L$. If  $X_1,...,X_N$ are independent copies of $X$, then for every $1 \leq k \leq N$ and $2 \leq q \leq p$,
$$
\|\bigl(\sum_{i \leq k} (X_i^*)^2\bigr)^{1/2}\|_{L_q} \leq cL(\sqrt{k\log(eN/k)} + \sqrt{q}).
$$
\end{Lemma}

\proof The proof follows from a comparison argument, showing that up to the $p$-th moment, the `worst case' is when $X$ is a gaussian variable.

Let $V_1,....,V_k$ be independent, nonnegative random variables  and set $V_1^\prime,....,V_k^\prime$ to be independent and nonnegative as well. Observe that if  $\|V_i\|_{L_q} \leq L\|V^\prime_i\|_{L_q}$ for every $1 \leq q \leq p$ and $1 \leq i \leq N$, then
\begin{equation} \label{eq:comp-lemma-1}
\|\sum_{i=1}^k V_i\|_{L_p} \leq L \|\sum_{i=1}^k V^\prime_i\|_{L_p}.
\end{equation}
Indeed, consider all the integer-valued vectors $\vec{\alpha}=(\alpha_1,...,\alpha_k)$, where $\alpha_i \geq 0$ and $\sum_{i=1}^k \alpha_i =p$. There are constants $c_{\vec{\alpha}}$ for which
$$
\|\sum_{i=1}^k V_i\|_{L_p}^p =\E \bigl(\sum_{i=1}^k V_i \bigr)^p=
\E \sum_{\vec{\alpha}} c_{\vec{\alpha}} \prod_{i=1}^k V_i^{\alpha_i}=\sum_{\vec{\alpha}} c_{\vec{\alpha}} \prod_{i=1}^k \E V_i^{\alpha_i},
$$
and an identical type of estimate holds for $(V_i^\prime)$. \eqref{eq:comp-lemma-1} follows if
$$
\prod_{i=1}^k \E V_i^{\alpha_i} \leq L^p \prod_{i=1}^k \E (V_i^\prime)^{\alpha_i},
$$
and the latter may be verified because $\|V_i\|_{L_q} \leq L \|V_i^\prime\|_{L_q}$ for $1 \leq q \leq p$.

\vskip0.4cm

Let $G=(g_i)_{i=1}^k$ be a vector whose coordinates are independent standard gaussian random variables. If $V_i = X_i^{2}$ and $V_i^\prime = c^2L^{2}g_i^{2}$, then by \eqref{eq:comp-lemma-1}, for every $1\leq q \leq p$,
$$
\|\sum_{i=1}^k X_i^{2}\|_{L_q} \leq c^2L^2 \|\sum_{i=1}^k g_i^2\|_{L_q}=c^2L^2 \left(\E\|G\|_{\ell_2^k}^{2q}\right)^{1/q}.
$$
It is standard to verify that
$$
\E \|G\|_{\ell_2^k}^{2q} \leq c^{2q} (\sqrt{k} + \sqrt{q})^{2q},
$$
and therefore,
$$
\|\sum_{i=1}^k X_i^2\|_{L_q} \lesssim L^2 \max\{k, q\}.
$$
By a binomial estimate,
\begin{align*}
& Pr\Bigl(\sum_{i \leq k} (X_i^*)^2 \geq t^2\Bigr) \leq \binom{N}{k} Pr\Bigl(\sum_{i \leq k} X_i^2 \geq t^2 \Bigr)
\\
\leq & \binom{N}{k} t^{-2q}\|\sum_{i \leq k} X_i^2\|_{L_q}^q \lesssim \left(\frac{eN}{k}\right)^k t^{-2q} \cdot L^{2q} (\max\{k,q\})^q,
\end{align*}
and if $q \geq k\log(eN/k)$ and $t=e u L \sqrt{q}$ for $u \geq 1$ then
\begin{equation} \label{eq:monotone-large-coordinates}
Pr\Bigl( \bigl(\sum_{i \leq k} (X_i^*)^2\bigr)^{1/2} \geq  e u L \sqrt{q} \Bigr) \leq u^{-2q}.
\end{equation}
Hence, setting $q = k\log(eN/k)$, tail integration implies that
$$
\|(\sum_{i \leq k} (X_i^*)^2)^{1/2}\|_{L_q} \lesssim L \sqrt{k \log(eN/k)},
$$
and if $q \geq k \log(eN/k)$, one has
$$
\|(\sum_{i \leq k} (X_i^*)^2)^{1/2}\|_{L_q} \lesssim L \sqrt{q},
$$
as claimed.
\endproof

The second preliminary result we require also follows from a straightforward binomial estimate:
\begin{Lemma} \label{lemma:moments-small-coordinates}
Assume that $\|X\|_{(p)} \leq L$ and let $X_1,...,X_N$ be independent copies of $X$. Consider $s \geq 1$, $1 \leq q \leq p$ and $1 \leq k \leq N$ that satisfies that $k \log(eN/k) \geq q$. Then
$$
\|\bigl(\sum_{i>k} (X_i^*)^{s}\bigr)^{1/s}\|_{L_q} \leq c(s)L N^{1/s},
$$
for a constant $c(s)$ that depends only on $s$.
\end{Lemma}

\proof Clearly, for every $1 \leq i \leq N$ and $2 \leq r \leq p$,
\begin{equation*}
Pr\left( X_i^* \geq t \right) \leq \binom{N}{i} Pr^i\left(X \geq t \right) \leq \binom{N}{i} \left(\frac{\|X\|_{L_r}^r}{t^r}\right)^i
\leq \left(\frac{eN}{i} \cdot \frac{L^r r^{r/2}}{t^r}\right)^i.
\end{equation*}
Hence, if $t=L\sqrt{r} \cdot e u$ for $u \geq 2$ and $r= 3\log(eN/i)$, then
\begin{equation} \label{eq:monotone-small-coordinates}
Pr \left(X_i^* \geq u \cdot eL \sqrt{3 \log(eN/i)} \right) \leq u^{-3 i\log(eN/i)}.
\end{equation}
Applying the union bound for every $i \geq k$, it follows that for $u \geq 4$, with probability at least $1-(u/2)^{-3 k \log(eN/k)}$,
\begin{equation} \label{eq:monotone-small-coordinates-union-bound}
X_i^* \leq u \cdot eL\sqrt{3\log(eN/i)}, \ \ {\rm for \ every \ } k \leq i \leq N.
\end{equation}
On that event
$$
\bigl(\sum_{i \geq k} (X_i^*)^s\bigr)^{1/s} \leq c(s)u L N^{1/s},
$$
and since $k \log(eN/k) \geq q$, tail integration shows that
$$
\|\bigl(\sum_{i \geq k} (X_i^*)^s\bigr)^{1/s}\|_{L_q} \leq c_1(s)L N^{1/s}.
$$
\endproof

\noindent{\bf Proof of Lemma \ref{lemma:(p)-estimates}.} Recall that $q_0=2+4\eta$, that $\xi \in L_{q_0}$ and that
$$
W=\left|\frac{1}{\sqrt{N}} \sum_{i=1}^N \eps_i \xi_i X_i\right| \IND_{{\cal A}_u}.
$$
Note that for every $(a_i)_{i=1}^N \in \R^N$ and any integer $0 \leq k \leq N$,
\begin{equation} \label{eq:bernoulli}
\|\sum_{i=1}^N \eps_i a_i\|_{L_q} \lesssim \sum_{i \leq k} a_i^* + \sqrt{q} \bigl(\sum_{i > k} (a_i^*)^2 \bigr)^{1/2}
\end{equation}
where the two extreme cases of $k=0$ and $k=N$ mean that one of the terms in \eqref{eq:bernoulli} is $0$.

Set $r=1+\eta$ and put $\theta=1/q_0$. Since $(\eps_i)_{i=1}^N$ are independent of $(X_i,\xi)_{i=1}^N$ and using the definition of the event ${\cal A}_u$,
\begin{align*}
 N^{q/2} \E W^q  = & N^{q/2} \E \IND_{{\cal A}_u} \E_\eps W^q \leq c^q \E \IND_{{\cal A}_u} \Bigl( \bigl(\sum_{i \leq k} \xi_i^* X_i^* \bigr)^q + q^{q/2} \bigl(\sum_{i > k} (\xi_i^*)^2 (X_i^*)^2 \bigr)^{q/2} \Bigl)
\\
\leq & c^q u^q \E_X\Bigl(  \bigl(\sum_{i \leq k} (N/i)^{\theta} X_i^*\bigr)^q + q^{q/2} \bigl(\sum_{i > k} (N/i)^{2\theta} (X_i^*)^2 \bigr)^{q/2} \Bigr).
\end{align*}

By the Cauchy-Schwarz inequality,
$$
\bigl(\sum_{i \leq k} (N/i)^\theta X_i^*\bigr)^q \leq \bigl(\sum_{i \leq k} (N/i)^{2\theta}\bigr)^{q/2} \cdot \bigl(\sum_{i \leq k} (X_i^*)^2\bigr)^{q/2},
$$
and
$$
\sum_{i \leq k} (N/i)^{2\theta} = \sum_{i \leq k} (N/i)^{1/1+2\eta} \leq \frac{c_1}{\eta}N^{1/(1+2\eta)} k^{2\eta/(1+2\eta)} \leq \frac{c_1}{\eta}N.
$$
Therefore,
$$
\E \bigl(\sum_{i \leq k} (N/i)^\theta X_i^*\bigr)^q \lesssim \eta^{-q/2} N^{q/2} \E \bigl(\sum_{i \leq k} (X_i^*)^2 \bigr)^{q/2} =(*).
$$
Also, by H\"{o}lder's inequality for $r=1+\eta$ and its conjugate index $r^\prime$,
$$
\bigl(\sum_{i > k} (N/i)^{2\theta} (X_i^*)^2 \bigr)^{q/2} \leq \bigl(\sum_{i \geq k} (N/i)^{2\theta r} \bigr)^{q/2r} \cdot \bigl(\sum_{i \geq k} (X_i^*)^{2r^\prime}\bigr)^{q/2r^\prime}
$$
and
$$
\sum_{i \geq k} (N/i)^{2\theta r}= \sum_{i \geq k} (N/i)^{(1+\eta)/(1+2\eta)} \leq \frac{c_1}{\eta} N.
$$
Hence,
$$
\E \bigl(\sum_{i > k} (N/i)^{2\theta} (X_i^*)^2 \bigr)^{q/2} \lesssim \eta^{-q/2r} N^{q/2r} \E \bigl(\sum_{i>k} (X_i^*)^{2r^\prime}\bigr)^{q/2r^\prime}=(**).
$$
Let $k \in \{0,...,N\}$ be the smallest that satisfies $k \log(eN/k) \geq q$ (and without loss of generality we will assume that such a $k$ exists; if it does not, the modifications to the proof are straightforward and are omitted).

Applying Lemma \ref{Lemma:com-with-gaussians} for that choice of $k$,
$$
(*) \leq c^{q}\eta^{-q/2} N^{q/2} \cdot L^q(\sqrt{k\log(eN/k)} + \sqrt{q})^q \leq c_1^q \eta^{-q/2} L^q N^{q/2} q^{q/2}.
$$
Turning to (**), set $s=2r^\prime \sim \max\{\eta^{-1},2\}$ and one has to control
$$
\E \bigl(\sum_{i>k} (X_i^*)^{s}\bigr)^{q/s}
$$
for the choice of $k$ as above. By Lemma \ref{lemma:moments-small-coordinates},
$$
\E \bigl(\sum_{i>k} (X_i^*)^{s}\bigr)^{q/s} \leq c^q(s)L^q N^{q/s}= c_1^q(\eta)L^q N^{q/2r^\prime}.
$$
Therefore,
$$
(**) \leq c^q(\eta) L^q N^{q/2r} \cdot N^{q/2r^\prime} = c^q(\eta) L^q N^{q/2}.
$$
Combining the two estimates,
$$
 N^{q/2} \E W^q \leq N^{q/2} \cdot c^q(\eta) L^q  q^{q/2},
$$
implying that $\|W\|_{L_q} \leq c(\eta) L$.
\endproof

\noindent{\bf Proof of Theorem \ref{thm:main-formulation}.} By Lemma \ref{lemma:(p)-estimates}, for every $1 \leq j \leq n$, $\|W_j\|_{(p)} \leq c(\eta)L$, and thus, by Lemma \ref{lemma:W-j-union-bound}, with probability at least $1-c_1t^{-2\beta}$,
$$
W_j^* \leq c(\eta) t L \sqrt{\beta \log(en/j)} \ \ \ {\rm for \ every \ } 1 \leq j \leq n.
$$
Moreover, $Pr({\cal A}_u) \geq 1-2/u^{q_0}$; therefore, with probability at least $1-c_1t^{-2\beta}-2u^{-q_0}$, for every $1 \leq j \leq n$,
$$
Z_j^* \leq c(\eta)  tu L \|\xi\|_{L_{q_0}} \sqrt{\beta \log(eN/j)}.
$$
Hence, on that event and because the norm $\sup_{v \in V} |\inr{v,\cdot}|$ is $K$ unconditional,
$$
\sup_{v \in V} |\inr{Z,v}| \leq K c(\eta) \sqrt{\beta} tu L \|\xi\|_{L_{q_0}} \sup_{v \in V}|\inr{Z_0,v}|,
$$
for a fixed vector $Z_0$ whose coordinates are $(\sqrt{\log(en/j)})_{j=1}^n$. Observe that $|\inr{Z_0,e_j}| \lesssim \E g_j^*$, and thus
$$
\sup_{v \in V} |\inr{Z_0,v}| \leq K \sup_{v \in V} |\sum_{i=1}^n v_i \E g_i^* |.
$$
Therefore, by Jensen's inequality, with probability at least $1-t^{-2\beta}-2u^{-q_0}$,
$$
\sup_{v \in V} \left|\frac{1}{\sqrt{N}}\sum_{i=1}^N \eps_i \xi_i x_i(j) \right|=\sup_{v \in V} |\inr{Z,v}| \lesssim_K c(\eta) \sqrt{\beta} tu L \|\xi\|_{L_{q_0}} \E \sup_{v \in V}|\inr{G,v}|.
$$
And, fixing $\beta$ and integrating the tails,
$$
\E \sup_{v \in V} \left|\frac{1}{\sqrt{N}}\sum_{i=1}^N \eps_i \xi_i x_i(j) \right| \lesssim_{K,\eta,L} \|\xi\|_{L_{q_0}} \ell_*(V),
$$
as claimed.
\endproof

\section{Applications in Sparse Recovery} \label{sec:sparse}
Spare recovery is a central topic in modern statistics and signal processing, though the problem we describe below is far from its most general form. Because a detailed  description of the subtleties of sparse recovery would be unreasonably lengthy, some statements may appear a little vague.  For more information on sparse recovery we refer the reader to the books \cite{MR2807761,MR2829871,MR3100033}, which are devoted to this topic.

\vskip0.4cm

The question in sparse recovery is to identify, or at least approximate, an unknown vector $v_0 \in \R^n$, and to do so using relatively few linear measurements. The measurements one is given are `noisy', of the form
$$
Y_i=\inr{v_0,X_i}-\xi_i \ \ {\rm for} \ 1 \leq i \leq N;
$$
$X_1,...,X_N$ are independent copies of a random, isotropic vector $X$ and $\xi_1,...,\xi_N$ are independent copies of a random variable $\xi$ that belongs to $L_q$ for some $q>2$.

The reason for the name ``sparse recovery" is that one assumes that $v_0$ is sparse: it is supported on at most $s$ coordinates, though the identity of the support itself is not known. Thus, one would like to use the given random data $(X_i,Y_i)_{i=1}^N$ and select $\hat{v}$ in a wise way, leading to a high probability estimate on the {\it error rate} $\|\hat{v}-v_0\|_{\ell_2^n}$ as a function of the number of measurements $N$ and of the `degree of sparsity' $s$.

\vskip0.4cm

In the simplest recovery problem, $\xi=0$ and the data is noise-free. Alternatively, one may assume that the $\xi_i$'s are independent of $X_1,...,X_N$, or, in a more general formulation, very little is assumed on the $\xi_i$'s.

The standard method of producing $\hat{v}$ in a noise-free problem and when $v_0$ is assumed to be sparse is the {\it basis pursuit algorithm}. The algorithm produces $\hat{v}$, which is the point with the smallest $\ell_1^n$ norm that satisfies $\inr{X_i,v_0}=\inr{X_i,v}$ for every $1 \leq i \leq N$.

It is well known \cite{MR2373017} that if $X$ is isotropic and $L$-subgaussian, $v_0$ is supported on at most $s$ coordinates and one is given
\begin{equation} \label{eq:optimal-gauss-est}
N = c(L)s\log\left(\frac{en}{s}\right)
\end{equation}
random measurements $(\inr{X_i,v_0})_{i=1}^N$, then with high probability, the basis pursuit algorithm has a unique solution and that solution is $v_0$.

Recently, it has been observed in \cite{LM_compressed} that the subgaussian assumption can be relaxed: the same number of measurements as in \eqref{eq:optimal-gauss-est} suffice for a unique solution if
$$
\max_{1 \leq j \leq n} \|\inr{X,e_j}\|_{(p)} \leq L \ \ {\rm for} \ \ p \sim \log n.
$$
And, the estimate of $p \sim \log n$ happens to be almost optimal. There is an example of an isotropic vector $X$ with iid coordinates for which
\begin{equation} \label{eq:BP-nec}
\max_{1 \leq j \leq n} \|\inr{X,e_j}\|_{(p)} \leq L \ \ {\rm for} \ \ p \sim (\log n)/(\log \log n)
\end{equation}
but still, with probability $1/2$ the basis pursuit algorithm does not recover even a $1$-sparse vector $v_0$ given the same number of random measurements as in \eqref{eq:optimal-gauss-est}.

\vskip0.4cm

Since `real world' data is not noise-free, some effort has been invested in producing analogs of the basis pursuit algorithm in a `noisy' setup. The most well known among these procedures is the LASSO (see, e.g. the books \cite{MR2807761,MR2829871} for more details) in which $\hat{v}$ is selected to be the minimizer in $\R^n$ of the functional
\begin{equation} \label{eq:LASSO}
v \to \frac{1}{N} \sum_{i=1}^N (\inr{v,X_i}-Y_i)^2 + \lambda\|v\|_{\ell_1^n},
\end{equation}
for a well-chosen of $\lambda$.

Following the introduction of the LASSO, there have been many variations on the same theme -- by changing the penalty $\| \ \|_{\ell_1^n}$ and replacing it with other norms. Until very recently, the behaviour of most of these procedures has been studied under very strong assumptions on $X$ and $\xi$ -- usually, that $X$ and $\xi$ are independent and gaussian, or at best, subgaussian.

\vskip0.4cm

One may show that Theorem \ref{thm:main-formulation} can be used to extend the estimates on $\|\hat{v}-v_0\|_{\ell_2^n}$ beyond the gaussian case thanks to two significant facts:
\begin{description}
\item{$\bullet$} The norms used in the LASSO and in many of its modifications happen to have a $1$-unconditional dual: for example, among these norms are weighted $\ell_1^n$ norms, mixtures of the $\ell_1^n$ and the $\ell_2^n$ norms, norms that are invariant under permutations, etc.
\item{$\bullet$} As noted in \cite{LM_reg_comp1}, if $\Psi$ is a norm, $B_\Psi$ is its  unit ball and $\hat{v}$ is the minimizer in $\R^n$ of the functional
\begin{equation} \label{eq:reg-procedures}
v \to \frac{1}{N} \sum_{i=1}^N (\inr{v,X_i}-Y_i)^2 + \lambda\Psi(v),
\end{equation}
then the key to controlling $\|\hat{v}-v\|_{\ell_2^n}$ is the behaviour of
\begin{equation} \label{eq:multi-in-sparse}
\sup_{v \in B_\Psi \cap r B_2^n} \left|\frac{1}{\sqrt{N}}\sum_{i=1}^N \xi_i \inr{X_i,v} - \E \xi \inr{X,v} \right|,
\end{equation}
\end{description}
which is precisely the type of process that Theorem \ref{thm:main-formulation} deals with.

It follows from Theorem \ref{thm:main-formulation} that if $\xi \in L_q$ for some $q>2$, the expectation of \eqref{eq:multi-in-sparse} is the same as if $\xi$ and $X$ were independent and gaussian. Thus, under those conditions, one can expect the `gaussian' error estimate in procedures like \eqref{eq:reg-procedures}. Moreover, because of \eqref{eq:BP-nec}, the condition that linear forms exhibit a subgaussian growth of moments up to $p \sim \log n$ is necessary, making the outcome of Theorem \ref{thm:main-formulation} optimal in this context.

\vskip0.4cm
The following is a simplified version of an application of Theorem \ref{thm:main-formulation}. We refer the reader to \cite{LM_reg_comp1} for its general formulation, as well as for other examples of a similar nature.

Let $X$ be an isotropic measure on $\R^n$ that satisfies $\max_{1 \leq j \leq n} \|\inr{X,e_j}\|_{(p)} \leq L$ for $p \leq c_0\log(n)$. Set $\xi \in L_q$ for $q >2$ that is mean-zero and independent of $X$ and put $Y=\inr{X,v_0}-\xi$.

Given an independent sample $(X_i,Y_i)_{i=1}^N$ selected according to $(X,Y)$, let $\hat{v}$ be the minimizer of the functional \eqref{eq:LASSO}.

\begin{Theorem} \label{thm:intro-LASSO-est}
Assume that $v_0$ is supported on at most $s$ coordinates and let $0<\delta<1$.
If $\lambda= c_1(L,\delta)\|\xi\|_{L_q} \sqrt{\log(en)/N}$, then with probability at least $1-\delta$, for every $1\leq p \leq 2$
\begin{equation*}
\norm{\hat v-v_0}_p \leq c_2(L,\delta)\|\xi\|_{L_q} s^{1/p}\sqrt{\frac{\log(ed)}{N}}.
\end{equation*}
\end{Theorem}
The proof of Theorem \ref{thm:intro-LASSO-est} follows by combining Theorem 3.2 from \cite{LM_reg_comp1} with Theorem \ref{thm:main-formulation}.

\begin{footnotesize}
\bibliographystyle{plain}
\bibliography{biblio}
\end{footnotesize}

\end{document}